\documentclass[10pt,twocolumn]{IEEEtran}

\usepackage{url,color}
\usepackage{subfigure}
\usepackage{amsfonts,mathrsfs}
\usepackage{amssymb,amsmath}
\usepackage{verbatim}
\usepackage{acronym}
\usepackage{graphicx}

\newcommand{\remove}[1]{}

\def\1{{\bf 1}}
\def\A{\mathscr{A}}

\def\O{{\Omega}}
\def\o{{\omega}}

\def\F{\mathscr{F}}

\def\b{\beta}

\def\e{\epsilon}

\def\a{\alpha}

\def\o{\omega}
\def\O{\Omega}

\def\R{\mathbb{R}}

\def\re{\mathbb{R}}
\def\sg{{\partial}}

\def\dist{{\rm dist}}

\newcommand{\EXP}[1]{\mathsf{E}\!\left[#1\right] }
\newcommand{\prob}[1]{\mathsf{Prob}\left\{#1\right\}}

\def\be{\begin{equation}}
\def\ee{\end{equation}}

\def\la{{\langle}}
\def\ra{{\rangle}}

\newtheorem{theorem}{Theorem}

\newtheorem{assumption}{Assumption}

\newtheorem{lemma}{Lemma}

\newtheorem{remark}[theorem]{Remark}

\makeindex 

\title{Minibatch stochastic subgradient-based projection algorithms for solving 
convex inequalities}

\author{Ion Necoara~\thanks{I. Necoara is with the Department of Automatic Control and Systems Engineering,  University Politehnica Bucharest, 060042 Bucharest, Romania. E-mail: \tt{ion.necoara@acse.pub.ro}.} and Angelia Nedi\'{c}~\thanks{A. Nedi\'{c}  is with the Electrical, Computer and Energy Engineering Department at Arizona State University, Tempe, AZ, USA. E-mail: {\tt angelia.nedich@asu.edu}.} 
\thanks{This work  is supported by the Executive Agency for Higher Education, Research and Innovation Funding (UEFISCDI), Romania,  PNIII-P4-PCE-2016-0731, project ScaleFreeNet, no. 39/2017. } }

\makeatother

\begin{document}
\maketitle

\begin{abstract}
This paper deals with the convex feasibility problem, where the feasible set is  given as the intersection of a (possibly infinite) number of closed convex sets. We assume that each set is specified algebraically as a convex inequality,  where the associated convex function is  general (possibly  non-differentiable).  For finding a point satisfying all the convex inequalities  we design and analyze   random projection algorithms using special subgradient iterations and extrapolated stepsizes. Moreover,  the iterate updates are performed based on  parallel  random observations of several constraint components. For these minibatch stochastic subgradient-based projection methods we prove sublinear convergence results and, under some linear regularity condition for the functional constraints, we prove linear convergence rates.   We  also derive conditions under which these  rates  depend explicitly on the minibatch size.   To the best of our knowledge, this work is the first deriving conditions  that show when   minibatch  stochastic subgradient-based projection updates  have a better  complexity than their single-sample variants.
\end{abstract}


\section{Introduction} 
\noindent Finding a point in the intersection of a collection of closed convex sets, that is the convex feasibility problem, represents a modeling paradigm for solving many engineering and physics problems, such as optimal control \cite{Calafiore01,PatNec:17}, robust control \cite{Alamo2009},  sensor networks \cite{BlaHer:06},  image recovery \cite{Combettes97}, data compression \cite{LieYan:05}, neural networks \cite{StaYan:98}, machine learning \cite{KunBac:17}.  Projection methods are very attractive in applications since they are able to handle problems of huge dimension and with a very large number of convex sets in the intersection. Projection methods were first used for solving systems of linear equalities  \cite{Kac:37} and linear inequalities \cite{Motzkin54}, and then extended to general convex feasibility problems, e.g. in \cite{Bregman67, Com:97, Gubin67, cdc2010, NecRic18}. For example, the alternating projection algorithm, which represents one of the first iterative algorithms for feasibility problems, rely at each iteration on orthogonal projections onto given individual sets taken in  a random, cyclic or greedy order \cite{Deutsch83, Deutch06, Halperin62, cdc2010, Polyak01}.  Otherwise, if the projection method uses, at the current iteration, an average of multiple projections of  the current iterate onto a subfamily of  sets, then it can be viewed as a minibatch projection algorithm \cite{Bauschke01, Bauschke03, NecRic18, Combettes97}.   
The convergence properties and even the inherent limitations of  projection methods have been intensely analyzed over the last decades, as it can be seen e.g.~in \cite{Bauschke01, Bauschke03, Com:97, Deutch06, NecRic18, cdc2010, Polyak01} and the references therein. 

\noindent \textit{Contributions}.  
In this paper we consider  convex feasibility problems with (possibly) infinite intersection of constraints. In contrast to the classical approach, where the constraints are usually represented as intersection of simple sets, which are easy to project onto, in this paper we consider that each constraint set is  given as the level set of a convex but not necessarily differentiable function.  For finding a point satisfying all convex inequalities   we propose   projection algorithms using the Polyak's subgradient update  (see~\cite{Polyak01}). Moreover, the iterate updates are performed based on  parallel  random observations of several constraint components and novel (adaptive) extrapolated  stepsize strategies.   For these  minibatch stochastic subgradient-based projection methods we derive sublinear convergence results and, under some linear regularity condition for the functional constraints, we prove linear convergence rates. We  also derive conditions under which these  rates  depend explicitly on the minibatch size (number of sets we project at each iteration).  From our best knowledge, this work is the first deriving theoretical conditions in terms of  the geometric properties of the functional constraints  that explain when   minibatch  stochastic subgradient-based projection updates  have a better  complexity than their non-minibatch  variants. More explicitly,  the convergence estimates for our  parallel projection algorithms depend on the key parameters $L$ or $L_N$,  defined in  \eqref{eq:Lnk}, which determines whether  minibatching helps ($L, L_N <1$) or not ($L=L_N=1$) and how much (the smaller $L$ or $L_N$, the better is the complexity).  Our algorithms are applicable to the situation where the whole constraint set of the problem is not known in advance, but it is rather learned in time through observations. Also, these algorithms are of interest for convex feasibility  problems where the constraints are known but their number  is either large or not finite.

\noindent  \textit{Content}. 
In Section~\ref{sec-problem_assum} we introduce our feasibilty problem  and derive some preliminary results. In Section~\ref{sec-method},  we present the Polyak's stochastic subgradient projection  method \cite{Polyak01} and derive its convergence rate under more general assumptions  than those in~\cite{Polyak01}. In Section~\ref{sec-parallel}  we consider minibatch variants with (adaptive) extrapolated stepsizes and derive the corresponding convergence rates depending on the minibatch size.  We provide some concluding remarks in~Section~\ref{sec-concl}.

\noindent \textit{Notation}. We will deal with a finite dimensional space $\R^n$,  where a vector is viewed as a column vector.  We use $\la x,y\ra$  to denote the inner product of two vectors $x,y \in \R^n$, and use $\|x\|$ to denote the standard Euclidean norm.  A vector $s_f$ is a subgradient of a convex function $f:\R^n\to\R$ at a point $x$ if  $f(x)+\la s_f,y-x\ra \le f(y)$ for all $y\in\R^n$. The set of all subgradients is denoted by $\sg f(x)$ and we use $f^+(x)=\max (0, f(x))$. We write $\dist(x,Y) =\min_{y\in Y}\|x-y\|$ ($\Pi_Y[x]$) for the distance (projection) from a vector $x$  to a closed convex set $Y$.  We  abbreviate {\it almost surely} by a.s.\ and {\it independent identically distributed}  by i.i.d.
 

\section{Problem, assumptions and preliminaries}
\label{sec-problem_assum}
\noindent  We consider the feasibility problem associated with a collection  $\{X_\o \}_{  \o\in\O }$ of convex closed sets in $\re^n$, where the index set $\O$ may be infinite.  We assume that each set $X_\o$ is specified as a convex inequality:
\begin{equation}
\label{eq-set}
X_\o=\{x  \in   \text{dom}(g_\o)   \mid g_\o(x)\le0\}\qquad\hbox{for all }\o\in\O,
\end{equation}
where each function $g_\o :\R^n \to  \bar \R$ is closed and convex,  with convex domain $\text{dom}(g_\o)$.  Some examples of the collection $\O$ include a finite set $\O=\{1,\ldots,p\}$ for some integer   $p\ge1$, an infinite countable set   $\O=\{1,2,\ldots\}$, or a countably infinite set  such as a closed interval  $\O=[a,b] \subset \R$. The convex feasibility problem of our interest is: 
 \begin{equation}
 \label{eq-problem}
 \hbox{find }x^*\in Y\cap (\cap_{\o\in \O} X_\o), 
 \end{equation}
 where $Y\subseteq\R^n$ is a (nonempty) 
 closed convex set and the set $X_\o$ has the functional representation~\eqref{eq-set}.
 The set $Y$ is assumed to have a simple structure for the projection operation 
 such as a halpfspace,  box or a ball. In the absence of such a constraint, we simply let $Y=\R^n$. The sets $X_\o$ are assumed to be complex for the projection operation.
Let $X^*$ denote the set of feasible points:
\[X^* \triangleq Y\cap \left(\cap_{\o\in\O} X_\o\right).\]
We assume that problem~\eqref{eq-problem} has a solution,
which is formalized as follows, together with assumptions on $Y$ and $g_\o$:

\begin{assumption} 
\label{assum-solution}
The set $Y\subseteq\R^n$ is closed, convex and simple for projection, while the function   $g_\o:\R^n \to  \bar \R $ is lower semicontinuous, convex and $Y \subseteq \text{rel int  dom} (g_\o)$  for each $\o\in\O$.  Moreover, the set  $X^*=Y\cap (\cap_{\o\in \O} X_\o)$ is nonempty.
\end{assumption}

\vspace{0.2cm}

\noindent A convex function that is lower semicontinuous has closed epigraph  \cite{Roc}. Hence, under Assumption~\ref{assum-solution}, each set $X_\o$ is nonempty, closed and convex and consequently the intersection set $X^*$ is also closed and convex.   We also assume that a subgradient $\partial g_\o(x)$ is available at any point $x \in  Y$ for all $\o \in \O$ (we mean an arbitrary subgradient if the set of them is not a singleton).  In our convergence analysis  we will work with the non-negative convex function $g_\o^+(x)=\max (0, g_\o(x))$. This function has the property that $g_\o^+(x) = 0$ if and only if $g_\o(x)  \leq 0$.  Furthermore,  we will use the following basic  result, which captures the change due to a feasibility step. The proof  can be found in~\cite{Pol:69}.

\begin{lemma} 
\cite{Pol:69}
\label{lemma-basiter}
Let $g: \R^n \to \bar{\R}$ be a convex function and $Y  \subseteq  \text{int dom} (g)$ be a nonempty  closed convex  set.  Let $x \in Y$ be arbitrary and consider Polyak's  subgradient step defined by:
\[z = \Pi_{Y} \left[ x - \b  \frac{g^+(x)}{\|d\|^2}\, d\right],\]
where $\beta>0$, and  $d \in\partial g^+(x)$ if $g^+(x)>0$ and otherwise $d \ne 0$ is arbitrary. Then, for any $y \in Y$ such that $g^+(y)=0$, the following inequality holds:
\[ \|z - y\|^2 \le  \|x - y\|^2 - \b(2-\beta)\,\frac{(g^+(x))^2}{\|d\|^2}.\]
\end{lemma}

\noindent Regarding the direction $d$ in Lemma~\ref{lemma-basiter},  let us note that  the function $g^+$ has a nonempty subdifferential set $\sg g^+(x)$ for all $x \in Y$,  since $g$, and consequently $g^+$,  are defined over  $Y  \subseteq  \text{int  dom}(g)$ \cite{Roc}. When $g(x)>0$,  then  $g^+(x)>0$ and  $0 \not \in \sg g^+(x)$ and, hence, the point $v  = x - \b  \frac{g^+(x)}{\|d\|^2} d $ is well defined. When $g(x)\le 0$, we have $g^+(x)=0$ and $v  = x - \b  \frac{g^+(x)}{\|d\|^2} d =x$ for any $d \ne 0$. Hence, in this case, regardless of the choice of the direction $d \ne 0$, we will always have $v=x$. We also use the assumption that the subgradients of $g_\o$ are bounded on the closed convex set~$Y$.

\begin{assumption}\label{asum-subgrads}
There exists a scalar $M_g>0$ such that
\[ \|s\|\le M_g \hbox{ for all $s\in\sg g_\o(y)$, all $y\in Y$ and all $\o\in\O$.} \]
\end{assumption}

\noindent We note that Assumption~\ref{asum-subgrads} always holds provided that  $Y$ is bounded and  $\text{dom} (g_\o) = \R^n$  for all $\o \in \O$ (recall that  for any convex lower semicontinuous function $g : \R^n \to \bar{\R}$ we have $\sg g(y)$  bounded for any $y \in \text{int} \; \text{dom} (g)$ \cite{Roc}).

 
\section{Stochastic subgradient-based  projection method}
\label{sec-method}
\noindent A simple algorithm for solving the convex feasibility problem~\eqref{eq-problem}  can be constructed using a randomly selected convex inequality at each iteration.
Specifically, by viewing the set $\O$ as the outcome space of a random variable $\o$ with a given distribution
$\pi$,
at each iteration $k$ we draw a random sample $\o_k$ according to $\pi$,
and we process the constraint $X_{\o_k}$. The constraint $X_{\o_k}$ is processed
by taking a step toward reducing the infeasibility of $X_{\o_k}$ at the current iterate and following 
by the projection step to remain in the set $Y$. Formally, the stochastic subgradient-based  projection (\texttt{SSP}) algorithm with a random constraint selection 
has the following update:

\begin{center}
\framebox{
\parbox{8 cm}{
\begin{center}
\textbf{ Algorithm \texttt{SSP}  }
\end{center}
Choose $x^0 \in Y$ and stepsize $\b \in (0, 2)$. For $k\geq 0$ do:
\begin{align}
&\text{Draw  sample} \;   \o_k \sim \pi   \; \text{and update:} \nonumber \\
&  x_{k+1}=\Pi_Y\left[x_k - \b\frac{g^+_{\o_k}(x_k)}{\|d_k\|^2}\,d_k\right],  
\label{eq-method}
\end{align}
}}
\end{center}
where  $d_k$ is a subgradient of $g^+_{\o_k}(x)$ at $x=x_k$ if $g_\o(x_k) >0$ (i.e.,  $d_k \in \sg g_{\o_{k}}^+ (x_{k}) = \sg g_{\o_{k}}(x_{k})$), and $d_k=d$ for some arbitrary $d\ne 0$, otherwise.  The initial point $x_0\in Y$ is assumed to be random (independent of $\o_k$).  Note that update rule \eqref{eq-method} can be viewed as a stochastic subgradient step with a special stepsize choice, as considered by Polyak  in \cite{Pol:69},  for solving the convex optimization problem having the  objective function expressed in terms of expectation:  
\begin{align}  
\label{eq:sref}
\min_{x \in Y} G(x) \quad \left(:= \EXP{  g_{\o}^+ (x)  } \right).  
\end{align}
Hence, we refer to the update rule \eqref{eq-method} of \texttt{SSP}  as~\textit{Polyak's subgradient} iteration. Note that the objective function $G(x) = \EXP{  g_{\o}^+ (x)  }$ in \eqref{eq:sref} can be viewed as a measure of infeasibility on $Y$ for the feasibility problem \eqref{eq-problem}, since   for any $x \in Y$ such that $G(x)>0$ implies that $x$ is infeasible, while   $G(x) = 0$ for any feasible point $x \in X^*$. However, $G(x)=0$ and $x \in Y$ does not imply that $x \in X^*$. Hence, the two problems \eqref{eq-problem}  and \eqref{eq:sref} are not equivalent, not without additional assumptions on the probability distribution $\pi$ and the functions $(g_\o)_{\o \in \O}$.  In Assumption \ref{asum-regularmod} below we provide a sufficient condition for the equivalence of the  problems \eqref{eq-problem}  and \eqref{eq:sref}, that is any $x \in Y$ satisfying $G(x) = 0$ is equivalent to $x \in X^*$.

\noindent   Note that under the assumption that the set  $X^*$ has a nonempty interior,  the iterative process ~\eqref{eq-method} with a stepsize $\b$ depending on the radius of a ball contained in $X^*$  has been proposed and studied by Polyak~\cite{Polyak01}.   In this section we analyze the convergence behavior of  algorithm \texttt{SSP} under more general stepsize rules  and less conservative assumptions on the sets $X_\o$ than those in~\cite{Polyak01}.  For this, we first  introduce the sigma-field $\F_{k}$ induced by the history of the method, i.e., by the realizations of the initial point $x_0\in Y$ and all the variables  $\o_t^i$ up to (including) iteration~$k$.  Specifically,
  \[ \F_{k} = \{x_0\} \cup  \left\{\o_{t} \mid \;  0 \le t \le k \right\}\quad\hbox{for }  k \ge 0. \]
For notational convenience, we   define $\F_{-1}=\{x_0\}$.  Then, we also impose the following \textit{linear regularity}  assumption.
  
\begin{assumption}\label{asum-regularmod}
  There exists a constant $c  \in (0, \infty)$ such that almost surely 
   for all $y\in Y$ and  all $k\ge 1$
  \[\dist^2(y,X^*) \le c \, \left( \EXP{g^+_{\o_{k}}(y) \mid \F_{k-1}} \right)^2.\]
\end{assumption}

\noindent  From Jensen's inequality applied to the convex function $\gamma(u) = u^2$ we always have the relation  $ \left( \EXP{g^+_{\o_{k}}(y) \mid \F_{k-1}} \right)^2  \leq  \EXP{(g^+_{\o_{k}}(y))^2\mid \F_{k-1}} $. Hence, Assumption~\ref{asum-regularmod} implies also:   
\begin{align}
\label{eq_impa3}
\dist^2(y,X^*)  \le c \, \EXP{(g^+_{\o_{k}}(y))^2\mid \F_{k-1}}. 
\end{align}     
Note that all our convergence results will also hold if we replace the condition from Assumption~\ref{asum-regularmod}  with \eqref{eq_impa3}. However, to establish a relation between the problems \eqref{eq-problem}  and \eqref{eq:sref} it is more natural to consider the condition from Assumption~\ref{asum-regularmod}  than \eqref{eq_impa3}. More precisely, under  Assumption~\ref{asum-regularmod} it is obvious that the two problems \eqref{eq-problem}  and \eqref{eq:sref} are equivalent, that is $x \in X^*$ if and only if  $x \in Y$ and $G(x) =  \EXP{  g_{\o}^+ (x)  } = 0$.  We summarize this discussion in the next lemma.

\begin{lemma}
Under  Assumption \ref{assum-solution} and  Assumption~\ref{asum-regularmod}  the feasibility problem \eqref{eq-problem}  is equivalent to  the stochastic optimization problem \eqref{eq:sref}. 
\end{lemma}

\noindent Assumption~\ref{asum-regularmod} summarizes all the  information we need regarding the distribution $\pi$ of the random variables  $\o_{k}$ and the initial point $x_0$. Under this assumption we will prove below that the sequence $x_k$ approaches $X^*$ at a linear rate. In the absence of this assumption, we will only prove that  our measure for infeasibility  $G(x) =  \EXP{  g_{\o}^+ (x)  } $ evaluated in a proper  average sequence will converge to zero at a sublinear rate.   Linear regularity assumption is standard in the convex feasibility literature \cite{Burke1993,Guler1992,cdc2010,Ned11,NecRic18}. It is always valid   provided that the interior of the intersection over the arbitrary  index set  $\A$  has an interior point \cite{Polyak01} (recall that the convergence analysis in \cite{Polyak01} has been given under this setting).  However, Assumption~\ref{asum-regularmod} holds for more general sets.  For example, when each set $X_\o$ is given by  a linear inequality $a_\o^T x + b_\o \leq 0$, one can verify that  the  intersection of these halfspaces  over any arbitrary  index set  $\A$ is linearly regular provided that the sequence $(a_\o)_{\o \in \A}$ is bounded, see  \cite{Burke1993,FerNec:19}. Hence,  Assumption~\ref{asum-regularmod} is also satisfied in this case and  $c$ is proportional to the  Hofman constant of the corresponding polyhedral set   \cite{NecRic18,Ned11}.  Furthermore, Assumption~\ref{asum-regularmod} holds under a strengthened Slater condition  for the collection of convex functional constraints $(X_\o)_{\o \in \A}$, such as the generalized Robinson condition, as detailed in Corollary 2 of~\cite{LewPan98}.  Next lemma derives a relation between $c$ and $M_g$:

\begin{lemma}
\label{lema:c}
Let Assumptions~\ref{asum-subgrads} and \ref{asum-regularmod} hold. Then, we have:
\[  c M_g^2 \geq 1. \]
\end{lemma}

\begin{IEEEproof}
Let $y \in Y$ be such that $y \not \in X^*$. Then, there exists $\bar{\o} \in \O$ such that the convex function $g_{\bar{\o}}$ satisfies $g_{\bar{\o}}(y) >0$.  Consequently, for any $s_g(y) \in \partial g_{\bar{\o}}(y)$  we also have $s_g(y) \in \partial g_{\bar{\o}}^+(y)$, and using convexity of $g_{\bar{\o}}^+$, we obtain:
\begin{align*}
0 & = g_{\bar{\o}}^+(\Pi_{X^*}[y]) \geq g_{\bar{\o}}^+(y) + \la s_g(y), \Pi_{X^*}[y] - y \ra\\ 
& \geq g_{\bar{\o}}^+(y) - M_g \|\Pi_{X^*}[y] -y\|,
\end{align*}
or equivalently
\[  g_{\bar{\o}}^+(y)  \leq  M_g \|\Pi_{X^*}[y] -y\|.  \]
On the other hand for those $\o \in \O$ for which $g_\o(y) =0$ we automatically have
\[  0 = g_\o^+(y)  \leq  M_g \|\Pi_{X^*}[y] -y\|.  \]
In conclusion, for any $\o \in \O$ and $y \in Y$ it holds:
\[  g_\o^+(y)  \leq  M_g \|\Pi_{X^*}[y] -y\|.  \]
Using the preceding  inequality and Assumption~\ref{asum-regularmod}, we get:
\begin{align*}
& \dist^2(y,X^*)  = \|\Pi_{X^*}[y] -y\| ^2 \leq c  \left( \EXP{ g^+_{\o_{k}}(y) \mid \F_{k-1}} \right)^2 \\
& \leq   c  \left( \EXP{ M_g  \|\Pi_{X^*}[y] -y\|   \mid \F_{k-1}} \right)^2 = c M_g^2 \ \dist^2(y,{X^*}),
\end{align*}
which proves our statement $c M_g^2 \geq 1$.  
\end{IEEEproof}


\subsection{Convergence analysis}
\label{subsec-method}
\noindent In this section we investigate the convergence behavior of the iteration ~\eqref{eq-method} of Algorithm \texttt{SSP}.  We first prove some descent  relation for the iteration~\eqref{eq-method} under a general probability  $\pi$.  

\begin{lemma}
\label{lem-descent-ssp}
Let Assumptions~\ref{assum-solution}--\ref{asum-subgrads} hold. Let also $x_{k+1}$ be obtained from update \eqref{eq-method} for some $x_k  \in Y$ and  $\b \in (0,2)$. Then, we have the following descent in expectation: 
\begin{align} 
\label{eq-mid}
& \EXP{\dist^2(x_{k+1},X^*)  \mid \F_{k-1}} \nonumber \\
& \le   \dist^2(x_{k},X^*) -   \frac{\b(2 -\beta)}{M_g^2}  \EXP{ (g_{\o_{k}}^+(x_k))^2  \mid \F_{k-1}}.
\end{align}
\end{lemma}
   
\begin{IEEEproof}      
From the definition of $x_{k+1}$ in~\eqref{eq-method} 
      and Lemma~\ref{lemma-basiter}, we obtain 
      for all $y\in X^*$ (for which we have $g_{\o_{k}}^+(y)=0$ 
      for any realization of $\o_{k}$)  and all $k\ge 0$
      \be\label{eq-one}
      \|x_{k+1} - y\|^2 \le  \|x_{k} - y\|^2
       -\b(2-\beta)\,\frac{(g_{\o_{k}}^+(x_{k}))^2}{\|d_k\|^2}.\ee  
 \noindent       By Assumption~\ref{asum-subgrads}, we have
        $\|d_k\| \le M_g$, implying for $y \in X^*$
       \[       \|x_{k+1} - y\|^2  \le  \|x_{k} - y\|^2
       -\frac{\b(2 - \beta)}{M_g^2}\,(g_{\o_{k}}^+(x_k))^2.\]
       By taking the minimum over $y\in X^*$ on 
       both sides of the preceding inequality, we have:
       \be
       \dist^2(x_{k+1}, X^*) \! \le \!  \dist^2(x_{k}, X^*)
       - \frac{\b(2 - \beta)}{M_g^2} (g_{\o_{k}}^+(x_k))^2. \nonumber 
       \ee
       Using the conditional expectation on $\F_{k-1}$   in the previous relation,   we obtain almost surely:
       \begin{align*} 
      & \EXP{\dist^2(x_{k+1},X^*)  \mid \F_{k-1}} \\
       & \le   \dist^2(x_{k},X^*) -   \frac{\b(2 -\beta)}{M_g^2}  \EXP{ (g_{\o_{k}}^+(x_k))^2  \mid \F_{k-1}},
       \end{align*}
which concludes our proof.        
\end{IEEEproof}  

\noindent When we consider general probability distributions $\pi$ and the descent Lemma \ref{lem-descent-ssp}, the following expected sublinear convergence  can be derived  for the measure of infeasibility  $G(x) =  \EXP{  g_{\o}^+ (x)  } $ evaluated in an   average sequence:   
\begin{theorem}
\label{lem-xiter-subl}
Let Assumptions~\ref{assum-solution}--\ref{asum-subgrads} hold and $x_{k}$ be the sequence generated by algorithm  \texttt{SSP} with $\b\in(0,2)$ and $x_0 \in Y$. Let also define the average sequence $\hat x_k = \frac{1}{k}\sum_{j=0}^{k-1} x_j$. If  $\EXP{\|x_0\|^2}<\infty$, then, almost surely, we have:
\[  \left( \EXP{G(\hat x_k)} \right)^2  \le  \frac{M_g^2 \EXP{ \dist^2(x_{0}, X^*)}}{  \b(2-\b)k }   \quad \forall k \ge 1. \]
\end{theorem}
   
\begin{IEEEproof}      
Using in the descent \eqref{eq-mid}  the  Jensen's inequality applied to the convex function $\gamma(u) = u^2$ and the definition of $G(x) =  \EXP{  g_{\o}^+ (x)  }$,      we obtain almost surely:
       \begin{align*} 
      & \EXP{\dist^2(x_{k+1},X^*)  \mid \F_{k-1}} \\
       & \le   \dist^2(x_{k},X^*) -   \frac{\b(2 -\beta)}{M_g^2}  \EXP{ (g_{\o_{k}}^+(x_k))^2  \mid \F_{k-1}} \\
       & \le   \dist^2(x_{k},X^*) -   \frac{\b(2 -\beta)}{M_g^2}  \left( \EXP{ g_{\o_{k}}^+(x_k)  \mid \F_{k-1}} \right)^2\\
       & =  \dist^2(x_{k},X^*) -   \frac{\b(2 -\beta)}{M_g^2}  \left( G(x_k) \right)^2,
       \end{align*} 
Taking now expectation over the entire history, we get the recurrence:
 \begin{align*} 
      & \EXP{\dist^2(x_{k+1},X^*)} \\
       & \le \EXP{  \dist^2(x_{k},X^*)} -     \frac{\b(2 -\beta)}{M_g^2} \EXP{ \left( G(x_k) \right)^2}.
       \end{align*}    
Adding the previous relation for $j=0$ to $k-1$ and using simple manipulations, we obtain:
 \begin{align*} 
& \frac{M_g^2  \EXP{\dist^2(x_{0},X^*)}}{\b(2 -\beta) k}   \geq \frac{1}{k} \sum_{j=0}^{k-1}   \EXP{ \left( G(x_j) \right)^2} \\ 
& \geq \frac{1}{k} \sum_{j=0}^{k-1}  \left( \EXP{  G(x_j) } \right)^2 \geq  \left( \frac{1}{k} \sum_{j=0}^{k-1}   \EXP{  G(x_j) } \right)^2 \\  
& =  \left(  \EXP{ \frac{1}{k} \sum_{j=0}^{k-1}    G(x_j) } \right)^2 \geq  \left(  \EXP{   G(\hat x_k) } \right)^2,
       \end{align*}   
where in the last inequality we used that $G$ is non-negative and convex. Moreover, for the term $\EXP{\dist^2(x_{0},X^*)}$ to be finite, it is sufficient to assume that  $\EXP{\|x_0\|^2}<\infty$.   This concludes our proof.            
 \end{IEEEproof}

\noindent When we consider  probability distributions $\pi$ satisfying the linear regularity condition (Assumptions~\ref{asum-regularmod}), then linear convergence for the expected distance  of the iterates of  \texttt{SSP} to $X^*$ can be derived. \noindent  Moreover, under  the condition $\EXP{\|x_0\|^2}<\infty$, the sequence $(\dist^2(x_k,X^*))_{k \geq 0}$ convereges to 0 almost surely. To show these statements, we make use of the following lemma, which can be found in~\cite{Polyak87}(Lemma 10).

\begin{lemma}[\cite{Polyak87}]\label{lem-polyak}
Let $\{v_k\}$ be a sequence of nonnegative random variables, with $\EXP{v_0}<\infty$ 
and satisfying a.s.\
\[\EXP{v_{k+1}\mid v_0,\ldots,v_k} \le(1-\a_k)v_k+\xi_k
\qquad\hbox{for all $k\ge0$}, \]
where $\a_k\in  (0, 1] $ and $\xi_k\ge0$ for all $k\ge0$, and 
\[\sum_{k=0}^\infty \a_k=\infty,\qquad
\sum_{k=0}^\infty \xi_k<\infty,\qquad
\lim_{k\to\infty}\frac{\xi_k}{\a_k}=0.\]
Then, $\lim_{k\to\infty} \EXP{v_k}=0$ and
almost surely $\lim_{k\to\infty}v_k=0$. Additionally, for any $\e>0$ and any $k>0$,
\[\prob{v_\ell \le\e,\ \forall\ell\ge k} \ge 1- \e^{-1}\,\left(\EXP{v_k}+\sum_{\ell=k}^\infty\xi_\ell\right).\]
\end{lemma}

\noindent Based on this result, we have the following theorem showing linear convergence of the distance of the iterates of \texttt{SSP} to $X^*$ in both expectation and probability. 
\begin{theorem}
\label{prop-sconverge}
Let Assumptions~\ref{assum-solution}--\ref{asum-regularmod} hold and $x_{k}$ be  the sequence generated by algorithm  \texttt{SSP} with $x_{0} \in Y$ and  $\b\in(0,2)$.  If  $\EXP{\|x_0\|^2}<\infty$, then, almost surely, we have:
\begin{align*}
\EXP{\dist^2(x_k, X^*) }  \le q^{k}\,\EXP{\dist^2(x_{0},X^*)}\  \forall k\ge1,  
\end{align*}
where    $q=1- \frac{\b(2-\beta)}{cM_g^2} \in [0, 1)$. Moreover, almost surely $\lim_{k\to\infty}\dist(x_k,X^*)=0$ and for any $\e>0$ and any $k \geq 1$,   we have in probability that
   \[ \prob{\dist^2(x_\ell,X^*) \le \e, \forall \ell\ge k} \ge 1- \frac{q^{k}}{\e}\,\EXP{\dist^2(x_0,X^*)}.\]
\end{theorem}

\begin{IEEEproof}
Using the conditional expectation on $\F_{k-1}$, 
       by Assumption~\ref{asum-regularmod}, we obtain almost surely:
       \[\EXP{(g_{\o_{k}}^+(x_k))^2 \mid \F_{k-1}} \ge \frac{1}{c} \dist(x_k,X^*).\]
       Therefore, using this inequality in~\eqref{eq-mid}, we obtain a.s.:
       \be \label{lem-xiter-sequential}
      \EXP{\dist^2(x_{k+1},X^*)  \mid \F_{k-1}}
       \le  q\, \dist^2(x_{k},X^*),  \ee 
       where
       $q=1- \frac{\b(2-\beta)}{cM_g^2}$.  Since $\b\in(0,2)$ the value of $\b(2-\b)$ lies in the interval $(0,1]$.      By   Lemma \ref{lema:c} we always have  $cM_g^2 \geq 1$. Hence,        it follows that $\b(2-\beta)/(cM_g^2) \in (0,1]$, implying that $q \in [0, 1)$. By taking now the total expectation in \eqref{lem-xiter-sequential} we get  $\EXP{\dist^2(x_k,X^*)} \le  q \EXP{ \dist^2(x_{k-1},X^*)}$. Using this relation recursively, we get linear convergence in expectation for  the distance of the iterates  to $X^*$:
\begin{align}
\label{eq-expseq}
\EXP{\dist^2(x_k, X^*) }  \le q^{k}\,\EXP{\dist^2(x_{0},X^*)}\  \forall k\ge1. 
\end{align}
Moreover, for the term $\EXP{\dist^2(x_{0},X^*)}$ to be finite, it is sufficient to assume that  $\EXP{\|x_0\|^2}<\infty$.   Furthermore, from \eqref{lem-xiter-sequential} we also see that the sequence $(\dist^2(x_k,X^*))_{k \geq 0}$
satisfies the conditions of Lemma~\ref{lem-polyak} with 
$v_k=\dist^2(x_k,X^*)$, $\a_k=1- q$ and $\xi_k=0$.
By Lemma~\ref{lem-polyak}, it follows that almost surely $\lim_{k\to\infty}\dist(x_k,X^*)=0$,
and that for any $\e>0$ and any $k>0$ we have
 \[\prob{\dist^2(x_\ell,X^*) \le\e, \forall \ell\ge k} \ge 1- \e^{-1}\,\EXP{\dist^2(x_k,X^*)}.\]
 By using relation~\eqref{eq-expseq} in the preceding inequality,
 we obtain the stated probability relation.
\end{IEEEproof}

\noindent When the set $Y$ is compact, since $x_0$ is random with realizations in $Y$,   the value $\EXP{\dist^2(x_0,X^*)}$ can be upper bounded by the diameter of the set $Y$, $\max_{x,y\in Y}\|x-y\|^2$,  which can be useful in applying the probability estimate of Theorem~\ref{prop-sconverge}. In this case, we get a lower bound on $\prob{\dist^2(x_\ell,X^*) \le\e, \forall \ell\ge k}$ with $\max_{x,y\in Y}\|x-y\|^2$ instead of $\EXP{\dist^2(x_0,X^*)}$.

\begin{remark}
\noindent In the  infeasible case, i.e. $X^* = \emptyset$, using a diminishing stepsize and the same arguments as in \cite{Polyak01},  we get a similar sublinear convergence rate as in  Theorem~\ref{lem-xiter-subl}. Hence, in the sequel we omit the analysis of this case.
\end{remark}


\subsection{Related work}
\noindent The paper most related to the results we derived in Section \ref{subsec-method} is  \cite{Polyak01}. In \cite{Polyak01} Polyak proves, that   under the assumption that the set  $X^*$ has a nonempty interior,  the iterative process ~\eqref{eq-method} with a stepsize $\b$ depending on the radius of a ball contained in $X^*$ has finite convergence.   On the other hand, Theorem~\ref{prop-sconverge} proves linear convergence of the distance of the iterates to $X^*$ under a more general linear regularity condition (Assumption \ref{asum-regularmod}) and a stepsize $\b$ which does not require knowledge of the set $X^*$. Note that our linear regularity condition covers the case   when $X^*$ has nonempty interior.  Furthermore, under a diminishing stepsize $\beta_k$, Polyak proves in  \cite{Polyak01} sublinear convergence  of the infeasibility measure $G$ in  an average sequence. Theorem~\ref{lem-xiter-subl} proves a similar result but for a constant stepsize $\b \in (0, 2)$. It is also important to note that our convergence analysis from Section \ref{subsec-method} is different from \cite{Polyak01}.  

\noindent When the sets $X_\o$ are easy for projections and $Y=\R^n$, by letting $g_{\o}(x)=\dist(x,X_\o) = \| x - \Pi_{X_\o}[x] \| $ and since $ x - \Pi_{X_\o}[x]/ \| x - \Pi_{X_\o}[x] \| \in \sg g_{\o}^+(x) $, the update \eqref{eq-method} reduces to the random projection iteration studied e.g. in~\cite{NecRic18,cdc2010}:
\[    x_{k+1} = x_k - \beta (x_k  - \Pi_{X_{\o_{k}}}[x_k]).   \]
Hence, our approach is more general since it allows to tackle also sets $X_\o$, described as the level set of a convex function $g_\o$, which are not easy for projection, but for which we can compute efficiently  a subgradient of $g_\o$.


\section{Minibatch stochastic subgradient-based  projection method}
\label{sec-parallel}
\noindent As noted in the previous section, the random update~\eqref{eq-method}, where $\{\o_k\}$ is  an i.i.d.\ sequence drawn according to some  distribution $\pi$ over $\O$,  can be interpreted as a stochastic approximation method for the convex problem  \eqref{eq:sref}, where all the functions $g^+_w(x)$ have a set of common minima $X^*$.  However,  distributed implementations of stochastic approximation methods  have become  recently   the \textit{de facto} architectural choice for large-scale stochastic problems.  Therefore, in what follows we will consider a minibatch variant of the update~\eqref{eq-method}, with a  probability distribution for the minibatch selection and (adaptive) extrapolated stepsize $\beta_k \geq 2$. 
It is expected that using minibatches of samples with a parallel batch processing and extrapolated stepsizes  would  be beneficial for a subgradient-based iterative process. 
Motivated by this idea, we consider a variant of Algorithm \texttt{SSP} with a minibatch of size $N$, i.e.,  having the current  iterate $x_{k}$, we sample $N$ constraints in parallel, and update as follows:

\begin{center}
\framebox{
\parbox{8 cm}{
\begin{center}
\textbf{ Algorithm \texttt{M-SSP}  }
\end{center}
Choose $x^0 \in Y$ and stepsizes $\b_k > 0$. For $k\geq 0$ do:
\begin{subequations}
\label{eq-method-parallel}
 \begin{align}
   &\text{Draw  sample} \;  J_k=( \o_k^1,\cdots, \o_k^N) \sim \mathsf{P}   \; \text{and update:} \nonumber \\
   & z_{k}^i    = x_{k} - \beta_k\, \frac{g^+_{\o_{k}^{i}}(x_{k})}{\|d_{k}^{i}\|^2}\, d_{k}^{i} 
    \quad \hbox{for } i=1,\ldots, N,\label{eq-zseqpar}\\
   & x_{k+1} = \Pi_Y[\bar z_k],\qquad \bar z_k=\frac{1}{N} \sum_{i=1}^N z_{k}^i \label{eq-xuppar},
  \end{align}
  \end{subequations}
}}
\end{center}
where $d_k^i \in \sg g_{\o_{k}^{i}}^+ (x_{k-1}) = \sg g_{\o_{k}^{i}} (x_{k-1})$ if $g_{\o_{k}^{i}}(x_{k})>0$,
  and $d_k^i=d$ for some arbitrary $d\ne 0$ otherwise. The initial point $x_0$ is 
  assumed to be random with outcomes in the set $Y$.  We  also need to redefine the sigma-field $\F_{k}$ induced by the history of the method, i.e., 
  by the realizations of the initial point $x_0\in Y$ and all the variables  $\o_t^i$ up to (including) iteration~$k$.  Specifically,
  \[ \F_{k} = \{x_0\} \cup  \left\{\o_{t}^{j} \mid \;  0 \le t \le k, \; 1 \le j\le N \right\}\quad\hbox{for }k\ge0. \]
  
\noindent 
{\it  We  will assume in the rest of the paper that Assumption~\ref{asum-regularmod} holds under this new sigma-field $\F_{k}$ for each $\o_k^i$ instead of $\o_k$}.  
The random $N$-tuple $J_k=(\o_k^1,\ldots,\o_k^N)$ generated according to the probability distribution $\mathsf{P}$ can be dependent
conditionally on $\F_{k-1}$. One choice is to draw $N$ independent samples $\o_k^1,\ldots,\o_k^N \sim \pi$.  When the index set $\O$ is finite, 
the indices $\o_k^i$ can be chosen randomly with or without replacement (e.g., given the realizations  $\o_k^1=j_1,\ldots,\o_k^{i-1}=j_{i-1}$, the index $\o_k^i$ is random with realizations in  
$\O\setminus\{j_1,\ldots,j_{i-1}\}$). Another possibility is to partition the set $\O$ into $N$ disjoint sets,   
$\cup_{i=1}^N\O_i=\O$, and select each $\o_k^i$ according to the uniform distribution over $\O_i$.

\noindent We also  need to  specify   how to choose the variable stepsize $\b_k$.  For this, let us define the following key  parameters that will play an important  role in  the way we define $\b_k$ and in the convergence analysis of \texttt{M-SSP}:
\begin{align}
\label{eq:Lnk}
& {\cal L}_N(x;J) =   \left\| \frac{1}{N}\sum_{i=1}^N\frac{g_{\o^{i}}^+(x)}{\|d^i\|^2}\, d^i\right\|^2 \Big{/} 
\left(\frac{1}{N} \sum_{i=1}^N\frac{(g_{\o^{i}}^+(x))^2}{\|d^i\|^2} \right),\nonumber  \\   
&  L_N^k = {\cal L}_N(x_k;J_k), \quad   L_N = \max_{x \in Y, J \sim \mathsf{P}} {\cal L}_N(x;J), \\
& L = \max_{x \in Y} \left\| \EXP{ \frac{g_{\o}^+(x)}{\|d_{\o} \|^2} \, d_{\o} } \right\|^2 \Big{/} 
\EXP{\frac{(g_{\o}^+(x))^2}{\|d_{\o}\|^2} },\nonumber 
\end{align}
where $J=(\o^1,\ldots,\o^N)$ and $d_{\o} \in  \sg g_{\o} (x)$ if $g_{\o}(x)>0$ and $d_{\o}=d$ for some arbitrary $d \ne 0$ otherwise. In the previous definitions of $ {\cal L}_N(x;J)$ and $L$ we use the convention that $0/0=0$.  By the convexity of the squared norm and Jensen inequality, we always have $L, {\cal L}_N(x;J)  \leq 1$ for all $x \in Y$  and  $J \sim \mathsf{P}$. Hence, we also have $L_N^k \leq L_N \leq  1$ for all $k \geq 0$.  However, there are convex functions $g_\o$ for which $L, L_N < 1$, as proved e.g., in the next lemma.  

\begin{lemma}
\label{lemma:polytope}
Let problem  \eqref{eq-problem} be described by $p$ linear inequalities, i.e. the functions $g_{\o}$ are given by:
\[  g_\o(x) = a_\o^T x + b_\o \leq 0  \quad  \forall  \o \in \Omega=\{1, 2, \ldots, p\},   \]
where $\|a_\o\|  =1$ for all~$\o$.  Define the matrix $A= [a_1 \cdots a_p]^T$ and for any $J =\{ \o^1 \cdots \o^N \} \subset \Omega$, sampled according to some probability $\mathsf{P}$, let $A_{J}$ be the submatrix of $A$ with the rows indexed in $J$.  If  the  submatrices  $A_{J}$ have at least rank two for all samples $J \sim \mathsf{P}$, then $L_N$ satisfies:
\begin{align}
\label{eq:LNlin}
L_N \leq \max_{ J \in 2^\O, |J| = N, J \sim \mathsf{P}} \frac{\lambda_{\max} (A_{J} A_{J}^T)}{N} <1.
\end{align}
In particular, if we consider uniform probability $\pi$ for sampling $\o$ and $A$ has at least rank two, then $L$  satisfies:
 \begin{align}
\label{eq:Llin}
L \leq   \frac{\lambda_{\max} (A A^T)}{p} <1.
\end{align}
\end{lemma}

\begin{IEEEproof}
 Let $x \in Y$ be fixed and $J =\{ \o^1 \cdots \o^N \} \subset \Omega$ be a sample of indexes selected  according to probability $\mathsf{P}$. Let us also define $J^+ = \{ \o \in J:  a_\o^T x + b_\o > 0  \}$.  In order to perform a nontrivial update in  \texttt{M-SSP} we must have $J^+ \not = \emptyset$.  Let $A_{J^+}$ be the submatrix of $A$ having the rows indexed  in the set $J^+$.  With these notations and assuming, without loss of generality,  that  $\|a_\o\|  =1$ for all $\o$, then  ${\cal L}_N(x;J)$ can be written explicitly as (recall that $|J^+| \geq 1$ in order to have a nontrivial update in \texttt{M-SSP}, otherwise ${\cal L}_N(x;J) =0$):
\begin{align*}
{\cal L}_N(x;J) &  
= \left\| \sum_{\o \in J^+} (a_\o x + b_\o) a_\o \right \|^2 \Big{/} \left(N \sum_{\o \in J^+} (a_\o x + b_\o)^2 \right) \\
& = \left\|   A_{J^+}^T ( A_{J^+} x + b_{J^+})  \right\|^2/ \left(N \| A_{J^+} x + b_{J^+} \|^2\right) \\
& \leq \frac{\lambda_{\max} (A_{J^+} A_{J^+}^T)}{N}   \leq \frac{\lambda_{\max} (A_{J} A_{J}^T)}{N} \\
& < \frac{\text{Trace}(A_{J} A_{J}^T)}{N} =  1  \quad \forall k,
\end{align*}
where the first inequality follows from the definition of the maximal eigenvalue $\lambda_{\max}$ of a matrix, the second inequality follows from $J^+ \subseteq J$ and  the eigenvalue interlacing  theorem, and the third  inequality holds strictly provided that the submatrix $A_{J}$ has at least rank two.   This concludes our first statement. For the second statement we first observe that if we choose $\o$ uniformly random, then  
\begin{align*}
& \left\| \EXP{ \frac{g_{\o}^+(x)}{\|d_{\o} \|^2} \, d_{\o} } \right\|^2 \Big{/} 
\EXP{\frac{(g_{\o}^+(x))^2}{\|d_{\o}\|^2} } \\ 
& =  \left\| \sum_{\o \in \Omega^+} (a_\o x + b_\o) a_\o \right \|^2 \Big{/} \left(p \sum_{\o \in \Omega^+} (a_\o x + b_\o)^2 \right)\\
& \leq \frac{\lambda_{\max} (A A^T)}{p},
\end{align*}
where  $\Omega^+ = \{ \o \in \Omega:  a_\o^T x + b_\o > 0  \}$ and we consider  nontrivial $x$'s satisfying $|\Omega^+| \geq 1$.  The rest follows using the  same reasoning as above. 
\end{IEEEproof}

\noindent Note that $L_N$ is an approximation of $L$ (empirical risk). Moreover,  $L_N^k$ is an online approximation of $L_N$. For particular sampling rules we can compute $L_N$ much more efficiently  than computing $L$, such as e.g., when we consider a uniform distribution over a fixed partition of $\Omega = \cup_{i=1}^\ell J_i$ of equal size sets. When $L_N$ is also difficult to compute we can use its online approximation $L_N^k$, whose computation is straightforward from the iteration of \texttt{M-SSP}.   Based on the  parameters $L$, $L_N$ and $L_N^k $  we define three strategies for the stepsize $\beta_k$: 
\begin{align*} 
&(\text{i}) \; \text{extrapolated stepsize} \;  \beta_k \in (0,  2/(1/N+(1-1/N)L)), \\  
&(\text{ii}) \; \text{minibatch extrapolated stepsize} \;  \beta_k \in (0,  2/L_N), \\  
&(\text{iii}) \;  \text{adaptive minibatch extrapolated stepsize} \; \beta_k \!\in\! (0, 2/L_N^k). 
\end{align*}   
From our best knowledge, these theree choices  for the stepsize in the minibatch  subgradient-based projection algorithm \texttt{M-SSP}  seem to be new.  
Moreover, since $L_N^k \leq L_N \leq 1$  it follows that $2/L_N^k \geq 2/L_N \geq 2$.  Similarly, since  $L \leq 1$  it follows that $1/N+(1-1/N)L \leq 1$ and thus $2/(1/N+(1-1/N)L) \geq 2$. However, when $L<1$ or $L_N <1$, we have $2/(1/N+(1-1/N)L) >1$ and  $2/L_N^k \geq 2/L_N > 2$, respectively.  Thus, we indeed   can choose  extrapolated stepsizes $\b_k > 2$ in the updates of \texttt{M-SSP}.   It is well-known that the  practical performance of projection methods  can be enhanced, and often dramatically so, using {\em extrapolation}, see e.g.,  \cite{Bauschke03,Com:97,NecRic18,Nec:19}.  In the next sections we also show  theoretically that our new extrapolated stepsizes bring benefits to the  algorithm \texttt{M-SSP} in terms of convergence~rates.


\subsection{Convergence analysis for extrapolated stepsize}
\noindent In this section we analyze the convergence behavior of algorithm \texttt{M-SSP} with  the extrapolated stepsize: 
 \[ \beta_k \in \left(0,  \frac{2}{ 1/N+(1-1/N)L } \right). \] First, we  prove some descent  relation for the iteration~\eqref{eq-method-parallel} of algorithm \texttt{M-SSP}. We consider  the  probability  $\mathsf{P}$  such that the $N$ samples  $\o_k^1,\ldots,\o_k^N$ are independent and drawn from the distribution  $\pi$.  
 
\begin{lemma}
\label{lem-xiter-parallel-L}
Let Assumptions~\ref{assum-solution}--\ref{asum-subgrads}  hold.  Let also $x_{k+1}$ be obtained from the update \eqref{eq-method-parallel} for some $x_k  \in Y$ and $\o_k^1,\ldots,\o_k^N$ are independent and drawn from the same distribution  $\pi$. Moreover,  let us also consider the  extrapolated stepsize  $\beta_k \in (0,  2/(1/N+(1-1/N)L))$.    Then, we have the following descent in expectation:
 \begin{align}
 \label{eq:lemaL}
    &\EXP{\dist^2(x_{k+1},X^*) \mid \F_{k-1}}  \le \dist^2(x_{k},X^*)  \\  
    & \quad  - \frac{\b_k(2N - \b_k ( 1 + (N-1)L))}{M_g^2 N} \EXP{ (g_{\o}^+(x_{k}))^2 \mid \F_{k-1}}. \nonumber 
    \end{align}   
   \end{lemma}
   
\begin{IEEEproof}
   By the projection non-expansiveness property, we have 
   $\|x_{k+1} - y\|^2\le \|\bar z_k - y\|^2$ for any $y\in X^*\subset Y$.
   Using this relation, we further get  that
\begin{align*}   
& \|x_{k+1} - y\|^2 \leq   \|  \frac{1}{N}\sum_{i=1}^N z_k^i - y\|^2 \\  
& = \left\|  x_k - y -    \frac{ \beta_k}{N}\sum_{i=1}^N   \frac{g^+_{\o_{k}^{i}}(x_{k})}{\|d_{k}^{i}\|^2}\, d_{k}^{i} \right \|^2\\
& =  \|  x_k - y\|^2 - 2  \frac{ \beta_k}{N} \sum_{i=1}^N    \frac{g^+_{\o_{k}^{i}}(x_{k})}{\|d_{k}^{i}\|^2} \langle  d_{k}^{i},   x_k - y \rangle \\ 
& \qquad + \beta_k^2  \left\|  \frac{ 1 }{N}\sum_{i=1}^N   \frac{g^+_{\o_{k}^{i}}(x_{k})}{\|d_{k}^{i}\|^2}\, d_{k}^{i}  \right \|^2.
\end{align*}
Now, using the convexity of $g^+_{\o_{k}^{i}}$  we get that
\[  0 =  g^+_{\o_{k}^{i}}(y) \geq g^+_{\o_{k}^{i}} (x^k) + \langle  d_{k}^{i},  y -   x_k  \rangle \quad \forall  y \in X^*,  \]
which, used in the previous derivations, yields:
\begin{align}
\label{eq:central}   
& \|x_{k+1} - y\|^2  \leq   \|  x_k - y\|^2 - 2  \frac{ \beta_k}{N} \sum_{i=1}^N    \frac{(g^+_{\o_{k}^{i}}(x_{k}))^2}{\|d_{k}^{i}\|^2}  \nonumber \\ 
& \qquad + \beta_k^2  \left\|  \frac{ 1 }{N}\sum_{i=1}^N   \frac{g^+_{\o_{k}^{i}}(x_{k})}{\|d_{k}^{i}\|^2}\, d_{k}^{i}  \right \|^2.
\end{align} 
From  \eqref{eq:central}  we further get:
\begin{align*}   
& \|x_{k+1} - y\|^2  \leq   \|  x_k - y\|^2 - 2  \frac{ \beta_k}{N} \sum_{i=1}^N    \frac{(g^+_{\o_{k}^{i}}(x_{k}))^2}{\|d_{k}^{i}\|^2} \\ 
& + \frac{\beta_k^2}{N^2}  \sum_{i=1}^N  \!  \frac{(g^+_{\o_{k}^{i}}(x_{k}))^2}{\|d_{k}^{i}\|^2} +    \frac{\beta_k^2}{N^2}  \sum_{i \not = j =1}^N \!\! \langle  \frac{g^+_{\o_{k}^{i}}(x_{k})}{\|d_{k}^{i}\|^2}\, d_{k}^{i} , \frac{g^+_{\o_{k}^{j}}(x_{k})}{\|d_{k}^{j}\|^2}\, d_{k}^{j}  \rangle   
\end{align*}
Minimizing both sides of the preceding inequality over $y\in X^*$, we find that
\begin{align*}
&\dist^2(x_{k+1},X^*)  \le  \dist^2(x_{k},X^*)  \\
& -  \b_k \left( 2-  \frac{\b_k}{N} \right)  \sum_{i=1}^N    \frac{1}{N}  \frac{(g^+_{\o_{k}^{i}}(x_{k}))^2}{\|d_{k}^{i}\|^2} \\
& +     \frac{\beta_k^2}{N^2}  \sum_{i \not = j =1}^N \!\! \langle  \frac{g^+_{\o_{k}^{i}}(x_{k})}{\|d_{k}^{i}\|^2}\, d_{k}^{i} , \frac{g^+_{\o_{k}^{j}}(x_{k})}{\|d_{k}^{j}\|^2}\, d_{k}^{j}  \rangle   
    \end{align*}
Taking the conditional expectation on $\F_{k-1}$ and using that $\o_k^1,\ldots,\o_k^N$ are independent and drawn from the same distribution  $\pi$,  we get for any $i \not = j$ that:
\begin{align*}
&\EXP{\dist^2(x_{k+1},X^*) \mid \F_{k-1}}  \le \dist^2(x_{k},X^*)  \\  
& -  \b_k \left( 2-  \frac{\b_k}{N} \right)    \EXP{ \frac{(g^+_{\o}(x_{k}))^2}{\|d_{k}^{\o}\|^2}  \mid \F_{k-1} } \\
& +     \frac{\beta_k^2}{N^2}  \sum_{i \not = j =1}^N \!\!  \langle  \EXP{ \frac{g^+_{\o}(x_{k})}{\|d_{k}^{\o}\|^2}\, d_{k}^{\o}  \mid \F_{k-1}} , \EXP{\frac{g^+_{\o}(x_{k})}{\|d_{k}^{\o}\|^2}\, d_{k}^{\o}  \mid \F_{k-1}}  \rangle, 
    \end{align*}   
where  for any $\o \in \Omega$ we define $d_k^{\o} \in  \sg g_{\o} (x_k)$ if $g_{\o}(x_k)>0$ and $d_k^{\o}=d$ for some arbitrary $d \ne 0$ otherwise.   Using now the definition of the constant $L$ we further  get:
\begin{align*}
&\EXP{\dist^2(x_{k+1},X^*) \mid \F_{k-1}}  \le \dist^2(x_{k},X^*)  \\  
& -  \b_k \left( 2-  \frac{\b_k}{N} \right)    \EXP{ \frac{(g^+_{\o}(x_{k}))^2}{\|d_{k}^{\o}\|^2}  \mid \F_{k-1} } \\
& +     \frac{\beta_k^2}{N^2} N(N-1)    \left\|  \EXP{ \frac{g^+_{\o}(x_{k})}{\|d_{k}^{\o}\|^2}\, d_{k}^{\o}  \mid \F_{k-1}} \right\|^2\\
&  \le \dist^2(x_{k},X^*)  \\ 
& -  \b_k \left( 2-  \frac{\b_k}{N}  -  \frac{\b_k}{N} (N-1)L\right)    \EXP{ \frac{(g^+_{\o}(x_{k}))^2}{\|d_{k}^{\o}\|^2}  \mid \F_{k-1} }. 
\end{align*}    
Note that  for $\beta_k \in (0,  2/(1/N+(1-1/N)L))$ the term  $2-  \b_k/N  -  (\b_k/N)(N-1)L \geq 0$. Hence, by combining the preceding recurrence with the assumption that the subgradients $d_k^\o$  are bounded (Assumption~\ref{asum-subgrads}),  we get our statement.
\end{IEEEproof}

\noindent From previous lemma we get the following sublinear convergence rate for the measure of infeasibility  $G(x) =  \EXP{  g_{\o}^+ (x)  } $ evaluated in an   average sequence:
\begin{theorem}
\label{mssp-xiter-subl-L}
Let assumptions of Lemma \ref{lem-xiter-parallel-L} hold with the constant extrapolated stepsize $\beta_k \equiv \beta = \frac{2 - \delta}{1/N+(1-1/N)L}$, where $\delta \in (0, 2)$, and $x_0 \in Y$, and define the average sequence $\hat x_k = \frac{1}{k}\sum_{j=0}^{k-1} x_j$. If  $\EXP{\|x_0\|^2}<\infty$, then, almost surely, we have for all $k \ge 1$:
\[  \left( \EXP{G(\hat x_k)} \right)^2  \le  \frac{ (1/N+(1-1/N)L)   M_g^2 \EXP{ \dist^2(x_{0}, X^*)}}{  \delta(2-\delta) k }. \]
\end{theorem}

\begin{IEEEproof} 
Since the samples   $\o_k^1,\ldots,\o_k^N$ are independent and drawn from the distribution  $\pi$, then we have: 
\begin{align*} 
& \EXP{  (g_{\o}^+(x_{k}))^2 \mid \F_{k-1}}  \geq  \left( \EXP{ g_{\o}^+(x_{k})  \mid \F_{k-1}} \right)^2 = (G(x_k))^2.  
\end{align*}
Using this relation in  \eqref{eq:lemaL} and  $\beta_k  =  \frac{2 - \delta}{1/N+(1-1/N)L}$, we get:
 \begin{align*}
    &\EXP{\dist^2(x_{k+1},X^*) \mid \F_{k-1}}  \le \dist^2(x_{k},X^*)  \\  
    & \quad  - \frac{\delta(2-\delta)}{(1/N+(1-1/N)L) M_g^2} (G(x_k))^2. 
    \end{align*}
Now, following the same reasoning as in the proof of Theorem \ref{lem-xiter-subl}, we get our statement.        
\end{IEEEproof} 

\noindent   When  additionally Assumptions~\ref{asum-regularmod} holds, then combining Lemma~\ref{lem-xiter-parallel-L} and Lemma~\ref{lem-polyak}, we obtain linear convergence rates in expectation and probability for  the expected distance  of the iterates of  \texttt{M-SSP} to $X^*$, with  the extrapolated stepsize.
   
\begin{theorem}
\label{prop-parconverge-L}
Let assumptions of Lemma \ref{lem-xiter-parallel-L} hold.  Let also  Assumption \ref{asum-regularmod} hold.  Also, assume that $c M_g^2 \geq  1/(1/N+(1-1/N)L)$ and $\EXP{\|x_0\|^2}<\infty$ for  $x_0 \in Y$. Then, the   sequence $(x_{k})_{k \geq 0}$ generated by the minibatch algorithm \texttt{M-SSP} with the constant extrapolated stepsize $\beta_k \equiv \beta = \frac{2 - \delta}{1/N+(1-1/N)L}$, where $\delta \in (0, 2)$, converges linearly in expectation:
  \begin{align*}
\EXP{\dist^2(x_k, X^*) }  \le q_{N,L}^{k}\,\EXP{\dist^2(x_{0},X^*)}\  \forall k\ge0,  
\end{align*} 
where $q_{N,L}=1-\frac{\delta(2-\delta)}{ (1/N+(1-1/N)L) cM_g^2} \in [0, 1)$,  and   almost surely    $\lim_{k\to\infty}\dist(x_k,X^*)=0$. Moreover, for any $\e>0$ and any $k>0$,   we also have the following convergence in probability: 
   \[ \prob{\dist^2(x_\ell,X^*) \le\e, \forall \ell\ge k} \ge 1- \frac{q_{N,L}^{k}}{\e}\,\EXP{\dist^2(x_0,X^*)}.  \]
\end{theorem}

\begin{IEEEproof}  Using the conditional expectation on $\F_{k-1}$,
       by Assumption~\ref{asum-regularmod}, we obtain almost surely:
       \begin{align*}   
       \EXP{  (g_{\o}^+(x_{k}))^2 \mid \F_{k-1}} & \geq  \left( \EXP{ g_{\o}^+(x_{k})  \mid \F_{k-1}} \right)^2 \\
       &  \ge \frac{1}{c} \dist(x_k,X^*).   
       \end{align*}
       Therefore, using this inequality in~\eqref{eq:lemaL}, we obtain a.s.:
    \begin{align*}
    &\EXP{\dist^2(x_{k+1},X^*) \mid \F_{k-1}}\\  
    & \le  \left( 1-\frac{\b_k(2N - \b_k( 1 + (N-1)L))}{cM_g^2 N} \right) \, \dist^2(x_{k},X^*).
    \end{align*}   
Using the expression of the  extrapolated stepsize $\beta_k = \frac{2 - \delta}{1/N+(1-1/N)L}$ in the previous relation, it follows that a.s.:   
   \be
   \label{eq-p1-L}
   \EXP{\dist^2(x_{k+1},X^*) \mid \F_{k-1}} \le  q_{N,L} \cdot \dist^2(x_{k},X^*),
   \ee
for all $k\ge0$,  with $q_{N,L}=1-\frac{\delta(2-\delta)}{ (1/N+(1-1/N)L) cM_g^2}$. Note that the conditions $cM_g^2  \geq 1/(1/N+(1-1/N)L)$, $1/N+(1-1/N)L \leq 1$ and $\delta \in (0, 2)$ implies that $q_{N,L} \in [0, 1)$.    Thus, the sequence $(\dist^2(x_k,X^*))_{k \geq 0}$ satisfies the conditions of Lemma~\ref{lem-polyak} with $v_k=\dist^2(x_k,X^*), \a_k=1- q_{N,L}$ and $\b_k=0$. Hence, it follows that $\lim_{k\to\infty}\EXP{\dist^2(x_k,X^*)}=0$  and that almost surely $\lim_{k\to\infty}\dist(x_k,X^*)=0$.  Also, for any $\e>0$ and any $k>0$ we have:
    \[\prob{\dist^2(x_\ell,X^*) \le\e, \forall \ell\ge k} \ge 1- \e^{-1}\,\EXP{\dist^2(x_k,X^*)}.\]
 By taking the total expectation in relation~\eqref{eq-p1-L}, we can see that
    $\EXP{\dist^2(x_k,X^*)} \le  q_{N,L} \EXP{ \dist^2(x_{k-1},X^*)},$
    which implies that for all $k\ge0$ we have a.s. linear convergence:
    \be\label{eq-p2-L}
    \EXP{\dist^2(x_k,X^*)} \le  q_{N,L}^k \EXP{\dist^2(x_{0},X^*)}.\ee
    Therefore, it also follows that 
     \[\prob{\dist^2(x_\ell,X^*) \le\e, \forall \ell\ge k} \ge 1- \frac{q_{N,L}^k}{\e}\,\EXP{\dist^2(x_0,X^*)},\]
     which concludes our statements.  
\end{IEEEproof}

\noindent  Regarding the assumption that $cM_g^2 \geq 1/(1/N+(1-1/N)L)$ in the previous theorem, we note that this assumption can be easily satisfied by choosing a larger value of $c$ or $M_g$ (since both of these values are defined in a form of upper bounds).    From  Theorems \ref{mssp-xiter-subl-L} and \ref{prop-parconverge-L}  we notice that  our convergence rates  depend on the minibatch size $N$  via the  term $1/N+(1-1/N)L$.   Note that if $L=1$, then $\beta = (2 - \delta)/(1/N+(1-1/N)L) = 2 - \delta \in (0, 2)$ and  $q_{N,L} = q$. Thus, in this case the convergence rates of \texttt{SSP} and  \texttt{M-SSP} are the same and   they do not depend on $N$. Hence, the complexity does not improve with minibatch size $N$.  However, as long as $L <1$ (and it can be also the case that $L \sim 0$),  then $q_{N,L}$ becomes  small, which shows that the minibatching algorithm \texttt{M-SSP}  with  \textit{extrapolated stepsize} has better performance than the  non-minibatch variant \texttt{SSP}.


\subsection{Convergence analysis for minibatch extrapolated stepsize}
\noindent In some cases we can easily compute $L$ (see e.g. Lemma \ref{lemma:polytope}). However, when it is difficult to compute $L$ we can use its empirical risk approximation $L_N$. Hence,   in this section we analyze the convergence behavior of algorithm \texttt{M-SSP} with  the minibatch extrapolated stepsize:  
\[ \beta_k \in \left(0,  \frac{2}{L_N} \right).  \] 
First, we  prove some descent  relation for the iteration~\eqref{eq-method-parallel} of algorithm \texttt{M-SSP} under a general probability  $\mathsf{P}$.  

\begin{lemma}
\label{lem-xiter-parallel}
Let Assumptions~\ref{assum-solution}--\ref{asum-subgrads}  hold.  Let also $x_{k+1}$ be obtained from the update \eqref{eq-method-parallel} for some $x_k  \in Y$ and for the  extrapolated stepsize  $\beta_k \in (0,  2/L_N)$.    Then, we have the following descent in expectation:
 \begin{align}
 \label{eq:lema6}
    &\EXP{\dist^2(x_{k+1},X^*) \mid \F_{k-1}}  \le \dist^2(x_{k},X^*)  \\  
    & \quad  - \frac{\b_k(2-\b_k L_N)}{M_g^2} \EXP{ \frac{1}{N}\sum_{i=1}^N(g_{\o_k^i}^+(x_{k}))^2 \mid \F_{k-1}}. \nonumber 
    \end{align}   
   \end{lemma}
   
\begin{IEEEproof}
Following the same proof as in Lemma \ref{lem-xiter-parallel-L} we get the following inequality (see  \eqref{eq:central}):  
\begin{align*}   
& \|x_{k+1} - y\|^2  \leq   \|  x_k - y\|^2 - 2  \frac{ \beta_k}{N} \sum_{i=1}^N    \frac{(g^+_{\o_{k}^{i}}(x_{k}))^2}{\|d_{k}^{i}\|^2}  \\ 
& \qquad + \beta_k^2  \left\|  \frac{ 1 }{N}\sum_{i=1}^N   \frac{g^+_{\o_{k}^{i}}(x_{k})}{\|d_{k}^{i}\|^2}\, d_{k}^{i}  \right \|^2,
\end{align*} 
for all  $y\in X^*$. Further, from the definition of $L_N^k$ and $L_N$ we have that:
\[   \left\|  \frac{ 1 }{N}\sum_{i=1}^N   \frac{g^+_{\o_{k}^{i}}(x_{k})}{\|d_{k}^{i}\|^2}\, d_{k}^{i}  \right \|^2 \leq L_N \left( \frac{1}{N} \sum_{i=1}^N    \frac{(g^+_{\o_{k}^{i}}(x_{k}))^2}{\|d_{k}^{i}\|^2} \right),  \]
which, used in the previous recurrence, yields:  
\begin{align*}   
& \|x_{k+1} - y\|^2  \leq   \|  x_k - y\|^2 - 2 \beta_k \left( \frac{ 1}{N} \sum_{i=1}^N    \frac{(g^+_{\o_{k}^{i}}(x_{k}))^2}{\|d_{k}^{i}\|^2} \right)  \\ 
& \qquad + \beta_k^2 L_N \left( \frac{1}{N} \sum_{i=1}^N    \frac{(g^+_{\o_{k}^{i}}(x_{k}))^2}{\|d_{k}^{i}\|^2} \right)  \\
& =   \|  x_k - y\|^2 - (2 \beta_k - \beta_k^2 L_N) \left( \frac{ 1}{N} \sum_{i=1}^N    \frac{(g^+_{\o_{k}^{i}}(x_{k}))^2}{\|d_{k}^{i}\|^2} \right).  
\end{align*} 
By combining the preceding recurrence with the assumption that the subgradients $d_k^i$  are bounded (Assumption~\ref{asum-subgrads}) and that $\beta_k \in (0,  2/L_N)$, we obtain  for all $y\in X^*$, 
\begin{align*}
   \|x_{k+1}-y\|^2 &\le  \|x_{k} - y\|^2  - \frac{(2 \beta_k - \beta_k^2 L_N)}{M_g^2}\,
    \frac{1}{N}\sum_{i=1}^N(g_{\o_k^i}^+(x_{k}))^2.
\end{align*}
    Minimizing both sides of the preceding inequality over $y\in X^*$, we find that
    \begin{align*}
    &\dist^2(x_{k+1},X^*) \\ 
    & \le  \dist^2(x_{k},X^*)  - \frac{\b_k(2-\b_k L_N)}{M_g^2} \frac{1}{N}\sum_{i=1}^N(g_{\o_k^i}^+(x_{k}))^2.
    \end{align*}
Taking the conditional expectation on $\F_{k-1}$,  we find:
    \begin{align*}
    &\EXP{\dist^2(x_{k+1},X^*) \mid \F_{k-1}}  \le \dist^2(x_{k},X^*)  \\  
    & \quad  - \frac{\b_k(2-\b_k L_N)}{M_g^2} \EXP{ \frac{1}{N}\sum_{i=1}^N(g_{\o_k^i}^+(x_{k}))^2 \mid \F_{k-1}},
    \end{align*}   
which concludes our statement.    
   \end{IEEEproof}

\noindent When we consider the probability distribution $\mathsf{P}$ such that the samples   $\o_k^1,\ldots,\o_k^N$ are independent and drawn from the distribution  $\pi$, the following expected sublinear convergence  can be derived  for the measure of infeasibility  $G(x) =  \EXP{  g_{\o}^+ (x)  } $ evaluated in an   average sequence:
\begin{theorem}
\label{mssp-xiter-subl}
Let Assumptions~\ref{assum-solution}--\ref{asum-subgrads} hold and the  samples   $\o_k^1,\ldots,\o_k^N$ are independent and drawn from the distribution  $\pi$. Let  also $x_{k}$ be the sequence generated by algorithm  \texttt{M-SSP} with the constant extrapolated stepsize $\beta_k \equiv \beta = \frac{2 - \delta}{L_N}$, where $\delta \in (0, 2)$, and $x_0 \in Y$, and define the average sequence $\hat x_k = \frac{1}{k}\sum_{j=0}^{k-1} x_j$. If  $\EXP{\|x_0\|^2}<\infty$, then, almost surely, we have:
\[  \left( \EXP{G(\hat x_k)} \right)^2  \le  \frac{L_N M_g^2 \EXP{ \dist^2(x_{0}, X^*)}}{  \delta(2-\delta) k }   \quad \forall k \ge 1. \]
\end{theorem}

\begin{IEEEproof} 
Since the samples   $\o_k^1,\ldots,\o_k^N$ are independent and drawn from the distribution  $\pi$, then we have: 
\begin{align*} 
& \EXP{ \frac{1}{N}\sum_{i=1}^N(g_{\o_k^i}^+(x_{k}))^2 \mid \F_{k-1}} \\ 
& \geq  \frac{1}{N}\sum_{i=1}^N \left( \EXP{ g_{\o_k^i}^+(x_{k})  \mid \F_{k-1}} \right)^2 = (G(x_k))^2.  
\end{align*}
Using this relation in  \eqref{eq:lema6} and  $\beta_k  = \frac{2 - \delta}{L_N}$, we get:
 \begin{align*}
    &\EXP{\dist^2(x_{k+1},X^*) \mid \F_{k-1}}  \le \dist^2(x_{k},X^*)  \\  
    & \quad  - \frac{\delta(2-\delta)}{L_N M_g^2} (G(x_k))^2. 
    \end{align*}
Now, following the same reasoning as in the proof of Theorem \ref{lem-xiter-subl}, we get our statement.        
\end{IEEEproof}

\noindent   When  we consider  probability distributions $\mathsf{P}$ satisfying only  the linear regularity condition (Assumptions~\ref{asum-regularmod}), i.e. there is no need to assume $\o_k^1,\ldots,\o_k^N$ to be independent, then combining Lemma~\ref{lem-xiter-parallel} and Lemma~\ref{lem-polyak}, we obtain linear convergence rates in expectation and probability for  the expected distance  of the iterates of  \texttt{M-SSP} to $X^*$, with  the extrapolated stepsize  $\beta_k \in (0,  2/L_N)$.
   
\begin{theorem}\label{prop-parconverge}
Let Assumptions~\ref{assum-solution}--\ref{asum-regularmod} hold.  Also, assume that $c M_g^2 \geq 1/L_N$ and $\EXP{\|x_0\|^2}<\infty$. Then, the   sequence $(x_{k})_{k \geq 0}$ generated by the minibatch algorithm \texttt{M-SSP} with the constant extrapolated stepsize $\beta_k \equiv \beta = \frac{2 - \delta}{L_N}$, where $\delta \in (0, 2)$, converges linearly in expectation:
  \begin{align*}
\EXP{\dist^2(x_k, X^*) }  \le q_N^{k}\,\EXP{\dist^2(x_{0},X^*)}\  \forall k\ge0,  
\end{align*} 
where $q_N=1- \frac{\delta(2-\delta)}{L_N cM_g^2} \in [0, 1)$,  and   almost surely    $\lim_{k\to\infty}\dist(x_k,X^*)=0$. Moreover, for any $\e>0$ and any $k>0$,   we also have the following convergence in probability 
   \[ \prob{\dist^2(x_\ell,X^*) \le\e, \forall \ell\ge k} \ge 1- \frac{q_N^{k}}{\e}\,\EXP{\dist^2(x_0,X^*)}.  \]
\end{theorem}

\begin{IEEEproof}  Using the conditional expectation on $\F_{k-1}$,
       by Assumption~\ref{asum-regularmod}, we obtain almost surely:
 \begin{align*}   
       \EXP{  (g_{\o_{k}^i}^+(x_{k}))^2 \mid \F_{k-1}} & \geq  \left( \EXP{ g_{\o_{k}^i}^+(x_{k})  \mid \F_{k-1}} \right)^2 \\
       &  \ge \frac{1}{c} \dist(x_k,X^*) \quad \forall i=1:N.   
       \end{align*}       
       Therefore, using this inequality in~\eqref{eq:lema6}, we obtain a.s.:
    \begin{align*}
    &\EXP{\dist^2(x_{k+1},X^*) \mid \F_{k-1}}\\  
    & \le  \left( 1-\frac{\b_k(2-\b_k L_N)}{cM_g^2} \right) \, \dist^2(x_{k},X^*).
    \end{align*}   
Using the expression of the  extrapolated stepsize $\beta_k \equiv \beta = \frac{2 - \delta}{L_N}$ in the previous relation, it follows that a.s.  for all $k\ge0$,
   \be
   \label{eq-p1}
   \EXP{\dist^2(x_{k+1},X^*) \mid \F_{k-1}} \le  q_N \cdot \dist^2(x_{k},X^*),
   \ee
with $q_N=1-\frac{\delta(2-\delta)}{L_NcM_g^2}$. The conditions $cM_g^2  \geq 1/L_N$, $L_N \leq 1$ and $\delta \in (0, 2)$ implies that $q_N \in [0, 1)$.    Thus, the sequence $(\dist^2(x_k,X^*))_{k \geq 0}$ satisfies the conditions of Lemma~\ref{lem-polyak} 
    with $v_k=\dist^2(x_k,X^*), \a_k=1- q_N$ and $\b_k=0$. Hence, it follows that $\lim_{k\to\infty}\EXP{\dist^2(x_k,X^*)}=0$  and that almost surely $\lim_{k\to\infty}\dist(x_k,X^*)=0$.  Also, for any $\e>0$ and any $k>0$ we have:
    \[\prob{\dist^2(x_\ell,X^*) \le\e, \forall \ell\ge k} \ge 1- \e^{-1}\,\EXP{\dist^2(x_k,X^*)}.\]
 By taking the total expectation in relation~\eqref{eq-p1}, we can see that
    $\EXP{\dist^2(x_k,X^*)} \le  q_N \EXP{ \dist^2(x_{k-1},X^*)},$
    which implies that for all $k\ge0$ we have linear convergence in expectation:
    \be\label{eq-p2}
    \EXP{\dist^2(x_k,X^*)} \le  q_N^k \cdot \EXP{\dist^2(x_{0},X^*)}.\ee
    Therefore, it also follows that 
     \[\prob{\dist^2(x_\ell,X^*) \le\e, \forall \ell\ge k} \ge 1- \frac{q_N^k}{\e}\,\EXP{\dist^2(x_0,X^*)},\]
     which concludes our statements.  
\end{IEEEproof}

\noindent  Regarding the assumption that $cM_g^2 \geq 1/L_N$ in the previous theorem, we note that this assumption can be easily satisfied by choosing a larger value of $c$ or $M_g$ (since both of these values are defined in a form of upper bounds).    From  Theorems \ref{mssp-xiter-subl} and \ref{prop-parconverge}  we notice that  our convergence rates  depend on the minibatch size $N$  via the key parameter $L_N$.   Note that if $L_N=1$, then $\beta = (2 - \delta)/L_N = 2 - \delta \in (0, 2)$ and  $q_N = q$. Thus, in this case the convergence rates of \texttt{SSP} and  \texttt{M-SSP} are the same and   they do not depend on $N$. Hence, the complexity does not improve with minibatch size $N$.  However, as long as $L_N <1$ (and it can be also the case that $L_N \sim 0$),  then $q_N$ becomes  small, which shows that the minibatching algorithm \texttt{M-SSP}  with \textit{minibatch extrapolated stepsize} has better performance than the  non-minibatch variant \texttt{SSP}.  
 

 \subsection{Convergence analysis for adaptive minibatch  extrapolated stepsize}
\noindent  If $L< 1$ or $L_N < 1$ and they  can be computed easily, then we have seen that  \texttt{M-SSP} algorithm with the (minibatch)  extrapolated steplengt  has (sub)linear convergence.  However, when $L$ or $L_N$ cannot be computed explicitly, we propose to approximate them online with $L_N^k$, i.e. we use  at each iteration an \textit{adaptive minibatch extrapolated  stepsize}  $\beta_k$ of the form:
\[  \beta_k \in \left( 0, \frac{2}{L_N^k}  \right),  \]
or equivalently, using the definition of $L_N^k$, as:
\begin{align*}
\beta_k  \in \left(0,      \frac{2}{N} \!\sum_{i=1}^N\!\frac{(g_{\o_{k}^{i}}^+(x_k))^2}{\|d_k^i\|^2} \Big{/} \left\| \frac{1}{N} \!\sum_{j=1}^N\!\frac{g_{\o_{k}^{j}}^+(x_k)}{\|d_k^j\|^2}\, d_k^j\right\|^2 \right),
\end{align*}
for any $x_k$ such that there exists at least one $i \in [1:N]$ satisfying $g_{\o_{k}^{i}}(x_k) >0$. Otherwise, we take $\beta_k \in (0,2)$. In this section we analyze the convergence behavior of algorithm \texttt{M-SSP} with  this adaptive choice for $\beta_k$. Note that the computational effort for computing  $L_N^k$ is the same as for the update in \eqref{eq-method-parallel}.  As in previous sections, we  first prove some descent  relation for the iteration~\eqref{eq-method-parallel} of algorithm \texttt{M-SSP} under a general probability  $\mathsf{P}$.

\begin{lemma}
\label{lem-xiter-parallel-adap}
Let Assumptions~\ref{assum-solution}--\ref{asum-subgrads} hold. Let also $x_{k+1}$ be obtained from the update \eqref{eq-method-parallel} for some $x_k  \in Y$ and for the adaptive minibatch extrapolated stepsize  $\beta_k =  \frac{2-\delta}{L_N^k} $ for some $\delta \in (0,2)$.    Then, we have the following descent in expectation: 
\begin{align} 
\label{eq:lema7}
&\EXP{\dist^2(x_{k+1},X^*) \mid \F_{k-1}}   \le   \dist^2(x_{k},X^*) \\ 
& \quad -  \frac{\delta (2 -  \delta)}{L_N  M_g^2} \EXP{ \frac{1}{N}\sum_{i=1}^N(g_{\o_k^i}^+(x_{k}))^2 \mid \F_{k-1}}.  \nonumber 
\end{align}
\end{lemma}
   
\begin{IEEEproof}
Following the same proof as in Lemma \ref{lem-xiter-parallel-L} we get \eqref{eq:central}, which we recall it here for convenience:   
\begin{align*}   
& \|x_{k+1} - y\|^2  \leq   \|  x_k - y\|^2 - 2  \frac{ \beta_k}{N} \sum_{i=1}^N    \frac{(g^+_{\o_{k}^{i}}(x_{k}))^2}{\|d_{k}^{i}\|^2}  \\ 
& \qquad + \beta_k^2  \left\|  \frac{ 1 }{N}\sum_{i=1}^N   \frac{g^+_{\o_{k}^{i}}(x_{k})}{\|d_{k}^{i}\|^2}\, d_{k}^{i}  \right \|^2,
\end{align*} 
for all  $y\in X^*$. Further, using the explicit expression for the adaptive extrapolated  stepsize $\b_k$, we obtain:
\begin{align*}   
& \|x_{k+1} - y\|^2  \leq   \|  x_k - y\|^2  - [2(2-\delta) - (2-\delta)^2] \\ 
& \left(\frac{1}{N} \sum_{i=1}^N    \frac{(g^+_{\o_{k}^{i}}(x_{k}))^2}{\|d_{k}^{i}\|^2} \right)^2   \left\|  \frac{ 1 }{N}\sum_{i=1}^N   \frac{g^+_{\o_{k}^{i}}(x_{k})}{\|d_{k}^{i}\|^2}\, d_{k}^{i}  \right \|^{-2}\\
&=   \|  x_k - y\|^2  - \frac{\delta (2-\delta)}{L_N^k}    \left(\frac{1}{N} \sum_{i=1}^N    \frac{(g^+_{\o_{k}^{i}}(x_{k}))^2}{\|d_{k}^{i}\|^2} \right). 
\end{align*} 
By combining the preceding recurrence with the assumption that the subgradients $d_k^i$  are bounded (Assumption~\ref{asum-subgrads}) and that $\delta \in (0,  2)$, we obtain  for all $y\in X^*$, 
\begin{align*}
   \|x_{k+1}-y\|^2 &\le  \|x_{k} - y\|^2  - \frac{\delta (2 -  \delta)}{L_N^k M_g^2}\,
   \left(  \frac{1}{N}\sum_{i=1}^N(g_{\o_k^i}^+(x_{k}))^2 \right).
\end{align*}
    Minimizing both sides of the preceding inequality over $y\in X^*$ and using that $L_N^k \leq L_N$ for all $k \geq 0$, we find that
    \begin{align*}
    &\dist^2(x_{k+1},X^*) \\ 
    & \le  \dist^2(x_{k},X^*)  - \frac{\delta (2 -  \delta)}{L_N M_g^2}\,
  \left(  \frac{1}{N}\sum_{i=1}^N(g_{\o_k^i}^+(x_{k}))^2 \right).
    \end{align*}
    Taking the conditional expectation on $\F_{k-1}$ in the previous relation we get our statement.    
   \end{IEEEproof}

\noindent When we consider the probability distribution $\mathsf{P}$ such that the samples   $\o_k^1,\ldots,\o_k^N$ are independent and drawn from the distribution  $\pi$, the following expected sublinear convergence  can be derived  for the measure of infeasibility  $G(x) =  \EXP{  g_{\o}^+ (x)  } $ evaluated in an   average sequence:

\begin{theorem}
\label{mssp-xiter-subl-adap}
Let Assumptions~\ref{assum-solution}--\ref{asum-subgrads} hold and the  samples   $\o_k^1,\ldots,\o_k^N$ are independent and drawn from the distribution  $\pi$. Let  also $x_{k}$ be the sequence generated by algorithm  \texttt{M-SSP} with the adaptive extrapolated stepsize  $\beta_k =  \frac{2-\delta}{L_N^k} $, where $\delta \in (0, 2)$, and $x_0 \in Y$, and define the average sequence $\hat x_k = \frac{1}{k}\sum_{j=0}^{k-1} x_j$. If  $\EXP{\|x_0\|^2}<\infty$, then  we have:
\[  \left( \EXP{G(\hat x_k)} \right)^2  \le  \frac{L_N M_g^2 \EXP{ \dist^2(x_{0}, X^*)}}{  \delta(2-\delta) k }   \quad \forall k \ge 1. \]
\end{theorem}

\begin{IEEEproof}
\noindent The proof of this theorem follows the same lines as in Theorem \ref{mssp-xiter-subl} and we omit it for brevity. 
\end{IEEEproof}

\noindent  Further,    when  we consider  probability distributions $\mathsf{P}$ satisfying the linear regularity condition (Assumptions~\ref{asum-regularmod}), i.e. there is no need to assume $\o_k^1,\ldots,\o_k^N$ to be independent, then combining Lemma~\ref{lem-xiter-parallel-adap} and Lemma~\ref{lem-polyak}, we obtain linear convergence rates in expectation and probability for  the expected distance  of the iterates of  \texttt{M-SSP} to $X^*$, with   the adaptive extrapolated stepsize  $\beta_k \in (0,  2/L_N^k)$.

\begin{theorem}\label{prop-parconverge-adap}
Let Assumptions~\ref{assum-solution}--\ref{asum-regularmod} hold.  Also, assume that $c M_g^2 \geq 1/L_N$ and $\EXP{\|x_0\|^2}<\infty$. Then, the   sequence $(x_{k})_{k \geq 0}$ generated by the minibatch algorithm \texttt{M-SSP} with the adaptive extrapolated stepsize $\beta_k  = \frac{2 - \delta}{L_N^k}$, where $\delta \in (0, 2)$, converges linearly in expectation:
  \begin{align*}
\EXP{\dist^2(x_k, X^*) }  \le q_N^{k}\,\EXP{\dist^2(x_{0},X^*)}\  \forall k\ge0,  
\end{align*} 
where $q_N=1- \frac{\delta(2-\delta)}{L_N cM_g^2} \in [0, 1)$,  and   almost surely    $\lim_{k\to\infty}\dist(x_k,X^*)=0$. Moreover, for any $\e>0$ and any $k>0$,   we also have the following convergence in probability 
   \[ \prob{\dist^2(x_\ell,X^*) \le\e, \forall \ell\ge k} \ge 1- \frac{q_N^{k}}{\e}\,\EXP{\dist^2(x_0,X^*)}.  \]
\end{theorem}

\begin{IEEEproof}
\noindent The proof  follows the same lines as in Theorem  \ref{prop-parconverge}. Hence, we omit it. 
\end{IEEEproof}


\subsection{When minibatching works?}
\noindent We notice, from  Theorems   \ref{mssp-xiter-subl-L}, \ref{mssp-xiter-subl} and \ref{mssp-xiter-subl-adap} on the one side  and  Theorems~\ref{prop-parconverge-L}, \ref{prop-parconverge} and \ref{prop-parconverge-adap} and the other side, that all three variants of  \texttt{M-SSP} using (adaptive minibatch) extrapolated stepsizes   have  (sub)linear convergence rates depending explicitly on the minibatch size $N$.  Moreover, when $L=1$ or  $L_N =1$,  the convergence rate of these variants of  \texttt{M-SSP} is the same as the one of the non-minibatch  method \texttt{SSP}, i.e of the form:
\begin{align*}
& \left( \EXP{G(\hat x_k)} \right)^2  \le  \frac{M_g^2 \EXP{ \dist^2(x_{0}, X^*)}}{  \b(2-\b) k }, \\
& \EXP{\dist^2(x_{k+1},X^*) \mid \F_{k-1}}  \le  q^k   \, \EXP{\dist^2(x_{0},X^*)},
\end{align*}
with $q=1-\frac{\b(2-\b)}{cM_g^2}$ and $\b \in (0,2)$. Hence, in this case, according to our results,  minibatching does not bring any benefits in terms of convergence rate. However, when $L<1$ or $L_N<1$ all  the three   variants of  \texttt{M-SSP} have (sub)linear convergence rates depending explicitly on  minibatch size $N$:  
\begin{align*}
& \left( \EXP{G(\hat x_k)} \right)^2  \le  \frac{{\cal L}_N \cdot M_g^2 \EXP{ \dist^2(x_{0}, X^*)}}{  \delta(2-\delta) k }, \\
& \EXP{\dist^2(x_{k+1},X^*) \mid \F_{k-1}}  \le  {\cal Q}_N^k   \cdot \EXP{\dist^2(x_{0},X^*)},
\end{align*}
where ${\cal L}_N$ is either $1/N + (1 - 1/N)L$ or $L_N$ and   ${\cal Q}_N = 1 - \frac{\delta(2 - \delta)}{{\cal L}_N c  M_g^2}$ and $\delta \in (0, 2)$.  For  example, for the linear rate (the analysis for the sublinear rate is similar)  if we choose the  optimal $\delta = 1$, we get ${\cal Q}_N = 1 -   \frac{1}{{\cal L}_N} \cdot \frac{1}{c M_g^2}$. Hence, ${\cal Q}_N$ is small provided that ${\cal L}_N \ll 1$. Furthermore, we observe that ${\cal Q}_N$ is with the order $1/{\cal L}_N$ better than $q$. In conclusion, as long as $L, L_N <1$ (and it can be also the case that $L, L_N \sim 0$),  then ${\cal Q}_N$ becomes  smaller than $q$, which shows that minibatching  improves complexity compared to single-sample variant.  Note that Lemma \ref{lemma:polytope} shows that e.g., polyhedral sets admit $L, L_N<1$.   To the best of our knowledge, this is the first time that Polyak's subgradient method with random minibatch   is shown to be better than its non-minibatch  variant.  We have identified $L$ and $L_N$ as the key quantities determining whether minibatching helps ($L, L_N < 1$) or not ($L, L_N = 1$), and how much (the smaller $L$ or $ L_N$, the more it helps).  

\noindent Note that  \texttt{M-SSP} algorithm does not require knowledge of   the subgradient norm  $M_g$, nor the constant $c$.  These values are only affecting the constants in the convergence rates, they are not needed for the stepsize selection.  Moreover, the adaptive minibatch  extrapolated stepsize $\beta_k = (2-\delta)/L_N^k$ can be easily implemented in practice even if the parameters $L, L_N$ are hard to compute.


\section{Conclusions}\label{sec-concl}
\noindent In this paper we have considered   
a convex feasibility problem with (possibly) infinite intersection of functional constraints.  For solving such a problem,  we  have proposed   
 minibatch stochastic subgradient methods motivated by Polyak's projection algorithm in~\cite{Polyak01}. At each iteration, our algorithms take  a subgradient step for minimizing the feasibility violation of the observed minibatch of constraints. The  updates are performed based on parallel  random observations of several  constraint components and based on (adaptive) extrapolated stepsizes. Under quite general conditions we have  derived sublinear rates, while  under  some additional  linear regularity condition on the functionals defining the sets, we have proved linear convergence rate  for this algorithm.   Moreover,  we have also derived conditions under which the  rate  depends explicitly on the minibatch size.  From our knowledge, this work is the first proving that  random minibatch  subgradient updates  have provably  better  complexity than their non-minibatch variants.

   
\bibliographystyle{amsplain}

\end{document}